\documentclass[11pt]{article}

\usepackage{amsthm,amssymb,latexsym,graphics,amscd,amsmath,amscd,enumerate,color}

\newcommand{\C}{\mathbb{C}}

\newcommand{\Cl}{\mathcal{C}\ell}

\newcommand{\mb}{\mathbb}
\newcommand{\mc}{\mathcal}
\newcommand{\R}{\mathbb{R}}

\newcommand{\Sc}{\mb{S}_{\mb{C}}}

\newcommand{\Spc}{\mbox{Spin}^{\mb{C}}}

\newcommand{\vol}{\operatorname{vol}}

\newtheorem*{cor}{Corollary}

\newtheorem{lem}{Lemma}

\newtheorem{prop}{Proposition}

\newtheorem*{mthm}{Main Theorem}

\begin{document}

\title{On the spectrum of the twisted Dolbeault Laplacian on line bundles over K\"ahler manifolds}
\author{Marcos Jardim \\ IMECC - UNICAMP \\
Departamento de Matem\'atica \\ Caixa Postal 6065 \\
13083-970 Campinas-SP, Brazil \\ and \\ Rafael F. Le\~ao \\ Departamento de Matem\'atica \\
Universidade Federal do Paran\'a \\ Caixa Postal 019081 \\ 81531-990 Curitiba-PR, Brazil}

\maketitle

\begin{abstract}
We use Dirac operator techniques to establish a sharp lower bound for the first eigenvalue of the twisted Dolbeault Laplacian on holomorphic line bundles over compact K\"ahler manifolds.

\end{abstract}


\section{Introduction}

Let $M$ be a compact K\"ahler manifold of complex dimension $n$, and let $E\to M$ be a holomorphic line bundle over $M$. Given a connection $\nabla_A$ on $E$, one has the decomposition:
\begin{equation}
  \nabla_A = \partial_A + \bar{\partial}_A ~~.
\end{equation}
The connection $\nabla_A$ is said to be compatible with the holomorphic structure on $E$ if $\bar{\partial}_A\bar{\partial}_A=0$; this is equivalent to requiring that the $(0,2)$-component fo the curvature 2-form $F_A$ vanishes.  The choice of a compatible connection $\nabla_A$ on $E$ induces a Hermitian metric on $E$.

If we denote the space of sections $\Gamma(E \otimes \wedge^{p,q}M)$ by $\Omega^{p,q}(E)$, then $\partial_A$ and $\bar{\partial}_A$ are first order differential operators acting as follows ($p,q=0,\dots,n$)
\begin{equation}
  \begin{split}
    \partial_A : \Omega^{p,q}(E) \rightarrow \Omega^{p+1,q}(E) \\
    \bar{\partial}_A : \Omega^{p,q}(E) \rightarrow \Omega^{p,q+1}(E)
  \end{split}
\end{equation}
Using the Hermitian structure induced on $E$ and the metric of $M$, we can define their formal adjoints
\begin{equation}
  \begin{split}
    \partial_A^* : \Omega^{p,q}(E) \rightarrow \Omega^{p-1,q}(E) \\
    \bar{\partial}_A^* : \Omega^{p,q}(E) \rightarrow \Omega^{p,q-1}(E)~.
  \end{split}
\end{equation}
These operators define a natural, second order differential operator on \linebreak $\Omega^{\bullet}(E) = \oplus_{p,q} \Omega^{p,q}(E)$, the so-called twisted Dolbeault Laplacian:
\begin{equation}
  \Delta_{\bar{\partial}_A} = \bar{\partial}_A \bar{\partial}_A^* + \bar{\partial}_A^* \bar{\partial}_A ~.
\end{equation}
Restricted to $\Omega^{0}(E)$, the twisted Dolbeault Laplacian simplifies to $\bar{\partial}_A^* \bar{\partial}_A$. This is the operator we shall concentrate on; from now on, we will denote it by $\Delta_A$. Notice that its kernel consists exactly of the holomorphic sections of $E$.

Finally, recall that a holomorphic line bundle $E$ is said to be negative if its first Chern class $c_1(E)$ can be represented as a closed $(1,1)$-form
$$ \frac{i}{2\pi} \sum_{k,l=1}^{n}
\varphi_{k\overline{l}}(p) dz^k\wedge d\overline{z}^l $$
such that the $n\times n$ matrix $\left[\varphi_{k\overline{l}}(p)\right]$ is negative definite for every $p\in M$ (see \cite{Ko}). Furthermore, the degree of $E\to M$ is defined as follows:
$$ \deg(E) = \int_M c_1(E)\wedge \omega^{n-1} ~~, $$
where $\omega$ is the K\"ahler form. It is not difficult to see that if $E$ is a negative line bundle, then $\deg(E)<0$.

Now, it follows from the Kodaira vanishing theorem is that if $E$ is a negative line bundle, then, for any compatible connection $\nabla_A$ on $E$, the twisted Dolbeault Laplacian $\Delta_A$ has trivial kernel (see \cite{Ko}).

Therefore, it is natural to pose the following problem. Let ${\cal A}(E)$ be the set of all compatible connections on a negative holomorphic line bundle $E\to M$. Consider the following functional
\begin{displaymath} \label{functional}
  \begin{split}
    \mathbf{\lambda} &: {\cal A}(E) \rightarrow \R^+ \\
    \nabla_A &\mapsto \lambda_1(\Delta_A)
  \end{split}
\end{displaymath}
that associates to each connection $\nabla_A\in{\cal A}(E)$ the value of the first nonzero eigenvalue of the associated twisted Dolbeault Laplacian $\Delta_A$. One would like to determine a lower bound for this functional, whether this lower bound is attained, and the characterization of its minima.

The goal of this paper is to provide the following solution to the problem just described. As usual, $\Lambda\alpha$ denotes contraction of a $(p,q)$-form $\alpha$ by the K\"ahler form of $M$. Recall that a connection $\nabla_A$ is said to be Hermitian-Einstein if the function $i\Lambda F_A$ is constant.

\begin{mthm}
Let $M$ be a compact K\"ahler manifold of complex dimension $n$, and let $E\to M$ be a holomorphic line bundle. If $\nabla_A$ is a compatible connection on $E$, then the eigenvalues of the twisted Dolbeault Laplacian $\Delta_A$ on sections of $E$ satisfies
\begin{equation}\label{lb}
  \lambda \geq - \frac{n}{2n-1}F_0
\end{equation}
where $F_0$ is the maximum of the function $i\Lambda F_A$.

Furthermore, if $\psi$ is an eigensection of $\Delta_A$ whose eigenvalue satisfies the equality, then $i\Lambda F_A$ is equal to the constant
\begin{equation}\label{c}
c=\frac{2\pi\deg(E)}{(n-1)!\vol(M)}
\end{equation}
and $\psi$ is in the kernel of the twistor operator. In particular, 
the connection $\nabla_A$ is the unique, up to gauge, Hermitian-Einstein connection on $E$.
\end{mthm}

Clearly, the estimate is meaningfull only when the function $i\Lambda F_A$ is strictly negative. For instance, if $E$ is a negative line bundle, then there exists a compatible connection $\nabla_A$ which satisfies this condition.

Notice that every holomorphic line bundle on a K\"hler manifold $M$ admits a compatible connection $\nabla_C$ which is unique up to gauge transformations and such that the function $i\Lambda F_C$ is constant \cite[p. 214]{DK}; such constant is fixed by Chern--Weil theory to be exactly $c$ as in (\ref{c}), see \cite{DK,Ko}.

In other words, the functional $\mathbf{\lambda}$ is bounded below by $-nc/(2n-1)$, and if this lower bound is attained by a connection $\nabla_A$, then $\nabla_A$ must be gauge equivalent to the Hermitian-Einstein connection on $E$. In particular, it follows from the discussion in \cite{JL2} that our estimate is sharp on every negative line bundle over a Riemann surface.

On the other hand it seems reasonable to ask whether the converse to the second part of our Main Theorem holds: if $\nabla_C$ is the Hermitian-Einstein connection on $E$, then the smallest eigenvalue of $\Delta_C$ is exactly $-nc/(2n-1)$. Based on previous work by Almorox--Prieto, the authors provided in \cite{JL2} a positive answer to this question for $M$ having dimension $1$, and one can construct examples of higher dimensional manifolds for which this converse statement is true. It would also be interesting to see whether this is also true for arbitrary higher dimensional K\"ahler manifolds.

Finally, our Theorem has an important consequence for the spectrum of the twisted complex Dirac operators on Riemann surfaces (case $n=1$), providing another example of the phenomenum discussed in \cite{JL3}.

\begin{cor}
Let $M$ be a Riemann surface, and let $E\to M$ be a holomorphic line bundle.
If $\nabla_A$ is a compatible connection on $E$, then the nonzero eigenvalues $\mu$ of the twisted complex Dirac operator $D_A$ satisfy
\begin{equation}\label{lb-d}
  \mu \geq \sqrt{-2F_0}~.
\end{equation}
where $F_0$ is the maximum of the function $i\Lambda F_A$.
\end{cor}

\paragraph{\bf Acknowledgments.}
The first named author is partially supported by the CNPq grant number 305464/2007-8 and the FAPESP grant number 2005/04558-0.


\section{Relation between the Dirac operator and the Dolbeault Laplacian}

The main component of the proof is the relation between the twisted complex Dirac operator and the twisted Dolbeault Laplacian on complex manifolds. Recall that every complex manifold is endowed with a canonical $\Spc$ structure associated to the complex structure. This structure comes with a complex Hermitian line bundle over $M$, called its determinant bundle. For the canonical $\Spc$ structure this bundle is just the canonical bundle of $M$, that is, $\wedge^{0,n}(M)$, where $n$ is the complex dimension of $M$. Furthermore, the spinor bundle associated to this canonical structure, denoted $\Sc$, can be identified with the bundle of holomorphic forms, i.e.
\begin{equation}
  \Sc \simeq \wedge^{0,*}(M),
\end{equation}
see \cite{N}. In this identification, the Clifford action on spinors can be explicitly given by
\begin{equation}
  \begin{split}
    \xi^k \cdot  \psi  = \sqrt{2} \xi^k \lrcorner \psi = \\
    \bar{\xi}^k \cdot \psi = \sqrt{2} \bar{\xi}^k \wedge \psi,
  \end{split}
\end{equation}
where $\psi \in \Gamma(\Sc)$, and $\{ \xi^k, \bar{\xi}^k \}$ is a unitary basis for $T^*M \otimes \C$.

If $E\to M$ is a holomorphic vector bundle, endowed with a compatible connection $\nabla_A$, we can consider the twisted spinor bundle $\Sc \otimes E$. Because $\Sc$ can be identified with $\wedge^{0,*}(M)$ the twisted spinor bundle can be identified with holomorphic forms with values in $E$, in other words,
\begin{equation}
  \Gamma \left( \Sc \otimes E \right) \simeq \Gamma \left( \wedge^{0,*}(M) \otimes E \right) = \Omega^{0,*}(E)
\end{equation}
The twisted spinor bundle has a natural structure of a Clifford module, with the Clifford action defined by
\begin{equation}
  \alpha \cdot \left( \psi \otimes s \right) = \left( \alpha \cdot \psi \right) \otimes s
\end{equation}
where $\psi \in \Gamma( \Sc )$, $s \in \Gamma( E )$ and $\alpha \in \Cl(T^*(M))$.

Using the connection $\nabla_A$on $E$, we can define the tensor product connection
\begin{equation}
  \nabla_{\tilde{A}} = \nabla_{\Sc} \otimes \mb{I}+\mb{I} \otimes \nabla_A
\end{equation}
on $\Sc \otimes E$. With this connection, and the natural module structure induced by the structure of $\Sc$, we define the twisted complex Dirac operator $D_A$ in the usual manner. With the identification between $\Gamma( \Sc \otimes E )$ and $\Omega^{0,*}(E)$, the complex twisted Dirac operator can be described in terms of the Cauchy operators of $E$, that is
\begin{equation}
  D_A = \sqrt{2} \left( \bar{\partial}_A + \bar{\partial}_A^* \right),
\end{equation}
which immediately gives the relation with the Dolbeault Laplacian
\begin{equation}
  \Delta_{\bar{\partial}_A} = \bar{\partial}_A \bar{\partial}_A^* + \bar{\partial}_A^* \bar{\partial}_A = \frac{1}{2} D_A^2.
\end{equation}
Here, $\Delta_{\bar{\partial}_A}$ and $D_A^2$ should be regarded as operators $\oplus_p\Omega^{0,p}(E)\to\oplus_p\Omega^{0,p}(E)$ which preserve the degree.


\section{Proof of the Main Theorem}

The idea of the proof is to use the Weitzenb\"ock formula for $D_A$ and the twistor equation to obtain the estimate. Recall that the Weitzenb\"ock formula for the twisted Dirac operator can be written as
\begin{equation}
  D_A^2 = \nabla_{\tilde{A}}^* \nabla_{\tilde{A}} + \frac{1}{4} R + \frac{1}{2} \Omega_{\Sc} + F_A
\end{equation}
where $\nabla_{\tilde{A}}^* \nabla_{\tilde{A}}$ is the trace Laplacian associated to the connection $\nabla_{\tilde{A}}$ on $\Sc \otimes E$, $R$ is the scalar curvature of $M$, $F_A$ is the curvature 2-form of $\nabla_A$ and $\Omega_{\Sc}$ is the curvature 2-form for some connection on the determinant bundle of the $\Spc$ structure (cf. \cite{LM}).

In principle, the connection on the determinant bundle of the $\Spc$ structure can be an arbitrary hermitian connection, but if we are dealing with the canonical $\Spc$ structure, associated to the complex structure of $M$, then the determinant bundle, in the above identification for $\Sc$, is $\wedge^{0,n}_M$ and the Chern connection of $M$ induces a canonical connection.

We are assuming that the connection $\nabla_A$ is compatible both with the Hermitian and holomorphic structures on $E$, and that the connection induced on $\wedge^{0,n}_M$ is constructed with these properties. So both curvature 2-forms, $\Omega_{\Sc}$ and $F_A$, are of type $(1,1)$. By \cite[Proposition 1]{JL2}, we know that the action of a $(1,1)$-form $\alpha$ on sections of $E$, which can be naturally identified with $\Omega^{0,0} E \subset \Sc \otimes E$, is given by
\begin{equation}
  \alpha \cdot \psi = -i ( \Lambda \alpha ) \psi,
\end{equation}
see \cite[Proposition 1]{JL2}.
Such characterization enables us to rewrite the Weitzenb\"ock formula for elements of $\Omega^{0,0}(E)$ as
\begin{equation}
  D_A^2|_{\Omega^0(E)} = 2 \bar{\partial}_A^* \bar{\partial}_A =
  \nabla_{\tilde{A}}^* \nabla_{\tilde{A}} + \frac{1}{4} R - 
  \frac{i}{2} \Lambda \Omega_{\Sc} - i \Lambda F_A
\end{equation}
To simplify the above formula, we need the folowing (see \cite[Proposition 2]{JL2}):

\begin{prop} \label{prop_cur_esc}
Let $M$ be a K\"ahler manifold and consider on the anti-canonical line bundle, $K_M^{-1} = \wedge^{0,n}_M$, the connection induced by the Chern connection of $M$. Let $\Omega_{\Sc}$ be the curvature 2-form of this connection, then we have
\begin{equation}
  i\Lambda \Omega_{\Sc} = \frac{R}{2}
\end{equation}
where $R$ is the Riemannian scalar curvature of $M$.
\end{prop}

Thus considering the connection on $\wedge^{0,n}_M$ induced by the Chern connection of $M$ reduces the above Weitzenb\"ock formula to
\begin{equation}
  D_A^2 = 2 \bar{\partial}_A^* \bar{\partial}_A = \nabla_{\tilde{A}}^* \nabla_{\tilde{A}} - i \Lambda F_A,
\end{equation}
which in turns can be simplified to
\begin{equation}\label{ww}
  \bar{\partial}_A^* \bar{\partial}_A = \frac{1}{2} \nabla_{\tilde{A}}^* \nabla_{\tilde{A}} - \frac{i}{2} \Lambda F_A.
\end{equation}


The twistor operator is given by
\begin{equation}
  \mc{T}_A = \sum_{k=1}^{m} e_k \otimes \left( \nabla_{\tilde{A},k} + \frac{1}{2n} e_k \cdot D_A \right) ~,
\end{equation}
where $\nabla_{\tilde{A},k}$ is the covariant derivative in the direction $e_k$ with respect to the product connection $\nabla_{\tilde{A}} = \nabla_{\Sc} \otimes \mb{I} + \mb{I} \otimes \nabla_A$ on $\Sc \otimes E$. It is a classical fact, see \cite{Ba}, that the twistor operator and the twisted Dirac operator can be related as
\begin{equation}
  \mc{T}_A^* \mc{T}_A = \nabla_{\tilde{A}}^* \nabla_{\tilde{A}} - \frac{1}{2n} D_A^2 ~.
\end{equation}

Now suppose that we have an eigensection $\psi \in \Omega^{0,0}(E)$ of $\bar{\partial}_A^* \bar{\partial}_A$, then using the Weitzenb\"ock formula (\ref{ww}) and taking the inner product with $\psi$ we obtain
\begin{equation}
  \lambda \mid \mid \psi \mid \mid^2 = \frac{1}{2} \mid \mid \nabla_{\tilde{A}} \psi \mid \mid^2 - \frac{1}{2} \int_M (i\Lambda F_A) \langle \psi \mid \psi \rangle
\end{equation}
Letting $F_0$ be the maximum value of the function $i\Lambda F_A$, the last integral can be estimated by
\begin{equation}\label{f0}
  -\frac{1}{2} \int_M (i\Lambda F_A) \mid \psi \mid^2 \geq -\frac{1}{2} F_0 \mid \mid \psi \mid \mid^2.
\end{equation}
Therefore we have
\begin{equation} \label{eq:frst_est}
  \lambda \mid \mid \psi \mid \mid^2 \geq \frac{1}{2} \mid \mid \nabla_{\tilde{A}} \psi \mid \mid^2 -\frac{1}{2} F_0 \mid \mid \psi \mid \mid^2.
\end{equation}

Now enters the relation between the twisted Dolbeault Laplacian and the twisted complex Dirac operator: if $\bar{\partial}_A^* \bar{\partial}_A \psi = \lambda \psi$ then we have $D_A^2 \psi = 2 \lambda \psi$, hence the equation relating the twistor operator and the Dirac operator, after taking the inner product with $\psi$, yields 
\begin{equation}\label{t1}
  \mid \mid \mc{T}_A \psi \mid \mid^2 = \mid \mid \nabla_{\tilde{A}} \psi \mid \mid^2 - \frac{\lambda}{n} \mid \mid \psi \mid \mid^2,  
\end{equation}
which implies that
\begin{equation}\label{t2}
  \mid \mid \nabla_{\tilde{A}} \psi \mid \mid^2 \geq \frac{1}{n} \lambda \mid \mid \psi \mid \mid^2.
\end{equation}
Using it in equation (\ref{eq:frst_est}) we obtain
\begin{equation}
  \lambda \geq - \frac{n}{2n-1} F_0 
\end{equation}
as desired.

Finally, remark that equality is attained if and only if the following two conditions hold:
\begin{itemize}
\item $i\Lambda F_A$ is constant, i.e. $\nabla_A$ is Hermitian-Einstein;
\item there exists an eigensection $\psi$ of the twisted Doulbeault Laplacian $\Delta_A$ such that $\psi\in\ker\mc{T}_A$.
\end{itemize}
The first condition is a consequence of equations (\ref{f0}) and (\ref{eq:frst_est}), while the second follows from (\ref{t1}) and (\ref{t2}).

The second fact requires a delicate analysis (compare with the untwisted case, \cite{F}). Notice however that $\psi\in\ker\mc{T}_A$ if and only if 
\begin{equation}
  \nabla_{\tilde{A}}^*\nabla_{\tilde{A}} \psi = \frac{1}{n} \lambda \psi
\end{equation}
i.e. the $\psi$ is a common eigensection of trace Laplacian $\nabla_{\tilde{A}}^*\nabla_{\tilde{A}}$ and Dolbeault Laplacian $\Delta_{A}$.

Finally, the Corollary follows easily from the Main Theorem and the following statement, proved in \cite{JL2}. Identifying $\Sc \otimes E \simeq \Omega^{0,*}(E)$, consider the projection operator
$p_0 : \Omega^{0,*}(E) \rightarrow \Omega^{0}(E)$.

\begin{lem}\label{dl}
If $\psi$ be an eigenstate of $D_A$, with non-null eigenvalue $\mu$, on a Riemann surface $\Sigma$, then we have
\begin{equation}
  p_0 \psi = \psi_0 \neq 0
\end{equation}
Furthermore, if $\mu$ is a nonzero eigenvalue of $D_A$, then $\frac{1}{2} \mu^2$ is an eigenvalue of $\bar{\partial}_A^* \bar{\partial}_A$.
\end{lem}
 

\end{document}